\documentclass[preprint,12pt]{elsarticle}
\usepackage{geometry}
\usepackage[normalem]{ulem}
\usepackage{comment}
\usepackage{graphicx,amsmath,amssymb,amsthm,marvosym,tikz}
\usepackage{nicematrix}
\usepackage{float}
\usepackage{multicol}
\usepackage[colorlinks]{hyperref}
\usetikzlibrary{patterns, shapes.misc, positioning}
\usetikzlibrary{decorations.pathreplacing, decorations.markings}
\newtheorem{theorem}{Theorem}[section]
\newtheorem{remark}{Remark}[section]
\newtheorem{lemma}{Lemma}[section]
\renewcommand{\d}{\ensuremath{\mathrm{d}}}
\newtheorem{corollary}{Corollary}[section]
\newtheorem{note}{Note}[section]

\newtheorem{definition}{Definition}[section]

\newtheorem{proposition}{Proposition}[section]
\journal{}
\begin{document}
\begin{frontmatter}
\title{On zero-divisor graph of the ring of Gaussian integers modulo $2^n$}
\author[a]{Aruna Venkatesan}
\ead{aruna_p220282ma@nitc.ac.in }	
\author[a]{ Krishnan Paramasivam}
\ead{sivam@nitc.ac.in}
\author[b]{M. Sabeel K}
\ead{sabeel.math@gmail.com}	
\address[a]{Department of Mathematics, National Institute of Technology Calicut \\ Kozhikode~\textnormal{673601}, India. }
\address[b]{Department of Mathematics, TKM College of Engineering  \\  Karicode, Kollam ~\textnormal{691005}, India.}

\begin{abstract}
For a commutative ring $R$, the zero-divisor graph, $\Gamma(R)$ is a simple graph having the vertex set as the set of all zero-divisors, $\mathcal{Z}(R)$ and two distinct vertices, $x$ and $y$ are adjacent if and only if $xy=0$. This article attempts to study the structure of the zero-divisor graph of Gaussian integers modulo $2^n$, focusing on its size, chromatic number, clique number, independence number, and matching through partitioning of zero-divisors. In addition, a few topological indices of the corresponding zero-divisor graph are found, such as the Wiener index, the Randi\'c index, the first Zagreb index, and the second Zagreb index.
\end{abstract}

\end{frontmatter}

\section{Introduction and Preliminaries}
In 1988, Beck \cite{MR0944156} first introduced the concept of the zero-divisor graph of a commutative ring. This concept was later redefined by Anderson and Livingston \cite{MR1700509}, who established the current definition. They defined the zero-divisor graph $\Gamma(R)$ of a commutative ring $R$ with vertices representing the non-zero zero-divisors $\mathcal{Z}(R)$, and two vertices $x$ and $y$ are adjacent if $xy=0$. In 2011, Osba et al. \cite{MR2458411} examined the zero-divisor graph of Gaussian integers modulo $n$, analyzing aspects such as the number of vertices, diameter, girth, and conditions under which the dominating number equals $1$ or $2$. They determined whether the zero-divisor graph of Gaussian integers is complete, complete bipartite, planar, regular, or Eulerian. Continuing their work, Osba et al. \cite{MR2783168} later explored conditions for $\Gamma(\mathbb{Z}_n[i])$ to be locally Hamiltonian or bipartite and determined its radius and chromatic number. Pirzada and Bhat \cite{article} investigated the clique number, connectivity, and degree conditions of $\Gamma(\mathbb{Z}_n[i])$. In 2022, Deepa and Kaur \cite {MR4370617} presented an algorithm for constructing the zero-divisor graph of the Gaussian integers modulo $2^n$ for $n 
\geq 1.$ They expressed $\Gamma(\mathbb{Z}_{2^n
}[i])$ as a generalized join graph  $G[G_1, G_2, \cdots, G_j]$ where each $G_j$ is either a complete graph (including loops) or its complement, and $G$ is the compressed zero-divisor graph of $\Gamma(\mathbb{Z}_{2^n}[i]).$\\
In this article, we study the structure of $\Gamma(\mathbb{Z}_{2^n}[i])$ in a different way from what is discussed in the literature. Element-wise analysis of $\Gamma(\mathbb{Z}_{n}[i])$ has been conducted in the literature except in \cite {MR4370617}, where the authors have given an algorithm to partition the zero-divisors of $\mathbb{Z}_{2^n}[i]$ and hence studied the corresponding zero-divisor graph. In this article, we conduct a substructure-wise analysis of the zero-divisor graph of $\mathbb{Z}_{2^n}[i]$ by partitioning the zero-divisors into associate classes.\\
First, we recall the definitions and notation used in the article. A graph $\Gamma=(V, E)$ is an ordered pair of two sets, where $V$ is a finite non-empty set of elements called vertices and $E$ is a subset of two elements of $V$ and is called the edge set. The cardinality of $V(\Gamma)$ is the order of $\Gamma$, and the cardinality of $E(\Gamma)$ is its size. A graph $\Gamma$ is connected if and only if a path exists between every pair of vertices $u$ and $v$. A graph on $n$ vertices in which any pair of distinct vertices is joined by an edge is a complete graph $K_n$. A complete subgraph of \( \Gamma \) with the largest order is called a maximal clique of \( \Gamma \), and the clique number \( \omega(\Gamma) \) is the number of vertices in such a maximal clique.
The number of edges incident on a vertex is called its degree, and a vertex of degree 1 is called a pendent vertex. The largest degree of a vertex is denoted by $\Delta$ and smallest degree is denoted by $\delta$. In a connected graph $\Gamma$, the distance between two vertices $u$ and $v$ of $\Gamma$ is the length of the shortest path between $u$ and $v$ in $\Gamma$. The connectivity $\kappa(\Gamma)$ of a graph $\Gamma$ is the smallest number of vertices whose removal from $\Gamma$ results in a disconnected graph or a trivial graph. In fact, for every graph $\Gamma$ of order $n$, $0 \leq \kappa(\Gamma) \leq n-1$. The eccentricity of a vertex $v$ of a connected graph $\Gamma$ is the distance from $v$ to the farthest vertex. The minimum eccentricity among the vertices of $\Gamma$ is the radius of $\Gamma$ and the maximum eccentricity is the diameter of $\Gamma$, denoted by $rad(\Gamma)$ and $diam(\Gamma)$ respectively.  A graph is planar if it can be redrawn in the plane so that no two edges cross. By a proper colouring of $\Gamma$, we mean an assignment of colours to the vertices of $\Gamma$, one colour to each vertex, so that adjacent vertices are coloured differently. The chromatic number $\chi(\Gamma)$ of $\Gamma$, is the minimum number of colours required to properly colour all the vertices of $\Gamma$. A graph $\Gamma$ is called perfect if every induced subgraph $\Gamma' \subseteq \Gamma$ has $\chi(\Gamma')=\omega(\Gamma')$. 
In 2006, Chudnovsky et. al. \cite{strongperfect} proved the characterization of perfect graphs.
A set of edges in a graph $\Gamma$ is independent if no two edges in the set are adjacent. By matching in $\Gamma$, we mean an independent set of edges in $\Gamma$. A matching $M$ saturates a vertex $u$, and $v$ is said to be $M$ - saturated if some edge of $M$ is incident with $v$; otherwise, $v$ is $M$-unsaturated. If every vertex of $\Gamma$ is $M$-saturated, the matching $M$ is perfect.
 $M$ is a maximum matching if $\Gamma$ has no matching $M'$ with $|M'|>|M|$ and the cardinality of such matching is denoted by $\alpha'(\Gamma).$ In particular, every perfect matching is a maximum.  An alternating path is a path whose edges are alternating between being in $M$ and not in $M$. An augmenting path $P$ with respect to a matching $M$ is an alternating path that starts and ends in unmatched vertices.
The smallest size of maximal matching is called the saturation number $s(\Gamma)$. 
Yannakakis and Gavril \cite{Yannakakis} proved that finding the smallest maximal matching is NP-hard even for bipartite (or planar) graphs with a maximum degree of 3. The saturation number is thus difficult to compute. However, it can be easily approximated by two factors.
 Every maximum matching is maximal, and therefore $s(\Gamma) \leq \alpha'(\Gamma)$, where $\alpha'(\Gamma)$ can be efficiently computed \cite{MR2536865}. Also see \cite{MR2074836}, \cite{MR3646699}, \cite{MR1981161} for more details on the bounds of $s(\Gamma)$. If $M$ is a maximal matching and $A$ is the set of end vertices of edges in $M$ then the set of vertices in $V(G)-A$ is an independent set of vertices in $\Gamma$ and hence $s(\Gamma) \geq \dfrac{|V(\Gamma)|- \alpha(\Gamma)}{2}$. By combining the above inequalities, $\dfrac{|V(\Gamma)|-\alpha(\Gamma)}{2} \leq s(\Gamma)\leq \alpha'(\Gamma).$
 One can refer to \cite{bondymurty} for more notation and terminology. We use $[j,k]$ to represent the set $\{j, j+1, j+2, \cdots ,k\}$, where $j<k$.

\begin{definition}\textnormal{\cite{cross}}
If $\alpha \in \mathbb{Z}[i]$, then the generalized Euler function from $\mathbb{Z}[i]$ to $\mathbb{N}, \phi_{\mathbb{Z}[i]}(\alpha)$, is defined to be the number of units in $\frac{\mathbb{Z}[i]}{\langle \alpha \rangle}$. That is, if $a+ib=up_1^{n_1}p_2^{n_2}\cdots p_k^{n_k}$, where $p_i$ is a prime in $\mathbb{Z}[i]$ and $u$ is a unit, then 
 $$\phi_{\mathbb{Z}[i]}(a+ib)=N(a+ib)\prod_{i=1}^k\bigl(1-\frac{1}{N(p_i)}\bigr).$$
 \end{definition}
 \begin{definition}\label{defn7} \textnormal{\cite{MR0130186}}
    Let $\Gamma$ be a given graph and $\{\Gamma_{\alpha}\}_{\alpha \in V(\Gamma)} $ be a collection of graphs indexed by $V(\Gamma)$. Then the generalized join of $\Gamma$ with $\{\Gamma_{\alpha}\}_{\alpha \in V(\Gamma)} $ is a graph $\tilde{\Gamma}$ with the vertex set $V(\tilde{\Gamma})=\{(x,y):x\in V(\Gamma) ~\textnormal{and}~ y\in V(\Gamma_x) \}$ and two vertices $(x,y)$ and $(x',y')$ are adjacent if and only if either $x$ is adjacent to $x'$ in $E(\Gamma)$ or $x=x'$ and $y$ is adjacent to $y'$ in $E(\Gamma_x)$.  If $\Gamma$ has $m$ vertices, then $\Gamma$ join of the collection $\{\Gamma_1,\Gamma_2, \ldots,\Gamma_m\}$ is denoted by $\Gamma[\Gamma_1,\Gamma_2,\ldots,\Gamma_m]$.
\end{definition}

\begin{definition}\textnormal{\cite{MR1700509}}\label{defn0}
Let $R$ be a commutative ring. A zero-divisor graph of the ring $R$ is a simple graph $\Gamma(R)$ having the vertex set as the set of all zero-divisors, $\mathcal{Z}(R)$ and two distinct elements, $x$ and $y$ of $\mathcal{Z}(R)$ are adjacent if and only if $xy=0$.
\end{definition}
 Let $R$ be a commutative ring with unity $1\neq0$, and let $x \in R$. The annihilator of $x$, denoted by $ann(x)$, is the set of elements of $R$ that annihilate $x$, that is, $ann(x) =\{y \in R : xy = 0\}$. We define a relation on $R$ by $x \sim y$ if and only if $ann(x) = ann(y)$. Clearly, this relation defines an equivalence relation on $R$ and therefore partitions $R$ into equivalence classes. By $[x]$, we shall denote the class of $x \in R$, that is, $[x]=\{y \in R : ann(x) = ann(y)\}$.
 
\begin{definition}\textnormal{\cite{spiroff}}\label{defn6}
For a commutative ring $R$ with $1 \neq 0$, a compressed zero-divisor graph of a
 ring $R$ is the undirected graph $\Gamma_E(R)$ with vertex set $R_E =\{[x]:x \in R\}$, where $[x]=\{y \in R : ann(x) = ann(y)\}$ and two distinct vertices
 $[x]$ and $[y]$ are adjacent if and only if $[x][y]=[0]=[xy]$, that is, if and only if $xy = 0$.
\end{definition}
For notation related to ring theory, we refer to \cite{gallian2021contemporary}.
Now we proceed to study a few topological indices, distance-based and degree-based graph indices. Wiener \cite{weiner} introduced the idea of the topological index while working on the boiling point of paraffin. The Wiener index of a connected graph $\Gamma$ is defined as the sum of distances between
each pair of vertices and is given by, 
\[
W(\Gamma)= \sum_{\alpha, \beta \in \Gamma} d(\alpha,\beta)
\]where $d(\alpha, \beta)$ is the length of the shortest path joining $\alpha$ and $\beta$. For more results and applications of the Wiener index of graphs, see \cite{gutman2,gutman3,hosaya}. The distance matrix $D(\Gamma)$ of a graph $\Gamma$ of order $n$, is an $n\times n$ matrix $(d_{jj'})$, where $d_{jj'}=d(\alpha_j,\alpha_{j'})$ for $j \neq j'$ and 0 otherwise. If $D(\Gamma)$ is the distance matrix of $\Gamma$, then the Wiener index of $\Gamma$ is given by,
\[W(\Gamma)=\frac{1}{2} \sum_{j=1}^{n} \sum_{j'=1}^{n}
d_{jj'}.
\]
The degree of vertex $\alpha$ of the zero-divisor graph $\Gamma$, denoted by $\text{d}(\alpha)$, is the number of vertices adjacent to $\alpha$.
The Randi\'c index (also known as the connectivity index) is a degree-based topological index. In 1976, Milan Randi\'c \cite{randic}
and is defined as,
\[
Rand(\Gamma)  =  \sum_{\alpha\beta \in E(\Gamma)} \frac{1}{(\text{d}(\alpha)\text{d}(\beta))^{1/2}}
\]
with summation over all pairs of adjacent vertices $\alpha$ and $\beta$ of the graph $\Gamma$. In 1972, Gutman and Trinajesti\'c \cite{gutman} introduced the Zagreb indices. For a graph $\Gamma$, the first Zagreb index $M_1(\Gamma)$ and the second Zagreb index $M_2(\Gamma)$ are, respectively, defined as follows:
\[
 M_1(\Gamma)=  \sum_{\alpha \in \mathcal{Z}(R)} (\text{d}(\alpha))^2
\]
\[
M_2(\Gamma)  =  \sum_{\alpha\beta \in E(\Gamma)} \text{d}(\alpha)\text{d}(\beta)
\]
For matrix-related notation, one can refer to \cite{bapat}.

\section{Zero-divisor graph of the ring of Gaussian integers modulo $2^n$}
The ring $\mathbb{Z}[i]=\{a+ib: a,b \in \mathbb{Z}\}$ of Gaussian integers is a  Euclidean domain with the Euclidean norm $N(a+ib)=a^2+b^2$ and hence is a principal ideal domain and unique factorization domain. An element $a + ib\in \mathbb{Z}[i]$ is prime in $\mathbb{Z}[i]$ if and only if $a^2 + b^2$ is prime in $\mathbb{Z}$. For an odd prime $p\in \mathbb{Z}$, can be expressed as $p=a^2 + b^2$ for some $a,b\in \mathbb{Z}$ if and only if  $p\equiv 1 \pmod{4}$ and  since for such $p, p=(a+ib)(a-ib)$, $p$ is not a prime in $\mathbb{Z}[i]$.  Therefore, an odd prime $p\in \mathbb{Z}$ is a prime in $\mathbb{Z}[i]$ if and only if $p\equiv 3 \pmod{4}$. In addition, $2$ is not a prime in $\mathbb{Z}[i]$, since $2$ can be factored into $(1+i)(1-i)$.
We consider the ring $\mathbb{Z}_n[i]$, which is isomorphic to the quotient ring of Gaussian integer modulo $n$ denoted by $\frac{ \mathbb{Z}[i]}{\langle n \rangle}$, that is, $$\mathbb{Z}_n[i]=\{ a+ib: a, b\in \mathbb{Z}_n \} \cong \frac{ \mathbb{Z}[i]}{\langle n \rangle}. $$

This section discusses important structural properties of the zero-divisor graph of the ring $R$, where $R$ is the ring of Gaussian integers modulo $2^n$. 
In this section, $R$ denotes the ring $\mathbb{Z}_{2^n}[i]$, where $n>1$. Since the prime factorization of $2^n=(-i)^n(1+i)^{2n}$, hence the proper divisors of $2^n$ in $\mathbb{Z}_{2^n}[i]$ are of the form $d_j=(1+i)^{2n-j}$ for $1\leq j \leq 2n-1.$ All these $d_{j}$'s are unique divisors up to associates. 

\begin{definition}\label{defn1}
Let $R$ be the ring of Gaussian integers modulo $2^n$, where $n>1$. For any proper divisor $d=a+ib$ of $2^n$ in $R$, we define  
$$V_{d}=\{ d' \in R : d'=ud \textnormal{~for some~} u \in \mathcal{U}(R)\}.$$
\end{definition}

Note that the above $V_{d}$ is a non-empty subset of $\mathbb{Z}_{2^n}[i]$ and $V_{d}$ is the associate class of $d$ in $\mathbb{Z}_{2^n}[i]$. In addition, the set of all $V_{d}$'s forms a partition of the set of all zero-divisors of $\mathbb{Z}_{2^n}[i]$. We now try to represent the structure of the zero-divisor graph of $R$ through these subsets $\{V_{d}: 1 <d< 2^n\} $.  Also note that for divisors $d$ and $d'$ of $2^n$ in $\mathbb{Z}_{2^n}[i]$, $V_{d}= V_{d'}$ if and only if $d$ and $d'$ are associates.

\begin{remark}{\label{note1}}
The prime factorization of $2^n$ in $\mathbb{Z}_{2^n}[i]$, is $2^n=(-i)^n(1+i)^{2n}$. Let $\{d_1, d_2, \ldots, d_{2n-1}\}$ denote the set of all proper divisors of $2^n$ in $\mathbb{Z}_{2^n}[i]$ with respect to this prime factorization. Then, $V_{d_1}, V_{d_2}, \ldots V_{d_{2n-1}}$ represent the distinct associate classes of the divisors $d_1, d_2, \ldots, d_{2n-1}$, respectively, in $\mathbb{Z}_{2^n}[i]$, where $d_j=(1+i)^{2n-j}$. 
\end{remark}

\begin{lemma}\label{lem1}
 Let $d_1, d_2, \ldots, d_{2n-1}$ be the set of distinct proper divisors up to associates of $2^n$ in the ring $\mathbb{Z}_{2^n}[i]$. Then,
 \begin{enumerate}
     \item[$(a)$] For any divisor $d_j$, the number of zero-divisors in  $V_{d_j}$ is $2^{j-1}$ .
     \item[$(b)$] A divisor $d_j$ is a nilpotent element of index $2$ if and only if $j \in [1,n]$.
     \item[$(c)$] The number of elements of nilpotency index $2$ is $2^n-1.$
     \item[$(d)$] If $d_j$ divides $d_{j'}$, then the number of zero-divisors in $V_{d_{j'}}$ divides the number of zero-divisors in $V_{d_j}$, where $j,j' \in [1,2n-1]$.
\item[$(e)$] The number of elements in $\cup_{j=1}^{2n-1} V_{d_j} $ is $2^{2n-1}-1.$
 \end{enumerate}
 \begin{proof}
     \begin{enumerate}
      \item[$(a)$] Clearly, $d_j=(1+i)^{2n-j}$. From \cite{young}, $|V_{d_j}|= \phi_{\mathbb{Z}[i]}(\frac{2^n}{d_j})=\phi_{\mathbb{Z}[i]}(1+i)^j= N(1+i)^{j-1}(N(1+i)-1)= 2^{j-1.}$
     \item[$(b)$] For any divisor $d_j$, we have 
 $d_{j}^2=  
((1+i)^{2n-j})^2 = 0$, if only if $j\in [1,n]$.
\item[$(c)$] Let $ud_j \in \displaystyle\cup_{j=1}^{n} V_{d_j}$. From $(b)$, $(u d_j)(u d_j)=(u)^2(d_j)^2=0$ implies that $ud_j$ is a nilpotent element of index 2. Using $(b)$, the elements from $\cup_{j=1}^{n} V_{d_j}$ are the only elements of index 2 in $\mathbb{Z}_{2^n}[i].$ Hence the number of elements of nilpotency index $2$ in $\mathbb{Z}_{2^n}[i]=\sum_{j=1}^{n} |V_{d_j}|=1+2+\cdots+2^{n-1}=2^n-1.$
\item[$(d)$] If $d_j$ divides $d_{j'}$ then $j < j'$. From $(a)$,  $|V_{d_j}|=2^{j-1}$ and $|V_{d_{j'}}|=2^{j'-1}$ and since $j < j'$, it is clear that $2^{j'-1}$ divides $2^{j-1}$. 
\item[$(e)$] $|\bigcup_{j=1}^{2n-1} V_{d_j}| = \sum_{j=1}^{2n-1} |V_{d_j}| = \sum_{j=1}^{2n-1} 2^{j-1} = 2^{2n-1}-1.$
     \end{enumerate}
 \end{proof}
\end{lemma}

\begin{lemma}\label{lem2}
An element $d_j$ is a nilpotent of index $2$ in $\mathbb{Z}_{2^n}[i]$ if and only if $\langle V_{d_j} \rangle$ is a clique in $\Gamma(\mathbb{Z}_{2^n}[i])$.
\end{lemma}
\begin{proof}
If $|V_{d_j}|=1$, then it is a clique. Let $|V_{d_j}| >1$. Assume $d_j$ is a nilpotent element of index 2 in $\mathbb{Z}_{2^n}[i]$.  Consider two elements $ud_j, u' d_j \in V_{d_j}$, where $u, u'$ are two units of $\mathbb{Z}_{2^n}[i].$ Then $ \langle V_{d_j} \rangle$ is a clique if and only if $uu' d_j^2=0$. Since $uu'\ne0$, $d_j^2=0$ is the only possibility.
\end{proof}
\begin{note} From Lemma $\ref{lem1}(b)$ and Lemma $\ref{lem2}$, it is clear that $\langle V_{d_j} \rangle$ is a clique in $\Gamma(\mathbb{Z}_{2^n}[i])$ if and only if $j \in [1,n].$ 
\end{note}
    
\begin{proposition} \label{prop1}
In the zero-divisor graph of $\mathbb{Z}_{2^n}[i]$, the subgraph induced by the set $V_{d_j}$ is the complete graph and the complement of the complete graph for $j \in [1,n]$ and $j \in [n+1,2n-1]$ respectively.     
\end{proposition}
\begin{proof}
   The first part of the Proposition directly follows from the previous Note. Let $u d_j, u' d_j$ be two elements of $V_{d_j}$ for $j \in [n+1,2n-1].$ Consider $(u d_j) \cdot (u' d_j)= (uu')(d_j)^2 \neq 0$ from Lemma $\ref{lem1}(b)$. Hence, for $j \in [n+1,2n-1]$, $\langle V_{d_j} \rangle$ is the complement of the complete graph.
   \end{proof}
\begin{corollary}\label{cor1}
    The clique number of each $\langle V_{d_j} \rangle$ is $2^{j-1}$ for $j \in [1,n]$.
\end{corollary}

\begin{proposition}\label{prop2}
The order of the zero-divisor graph of $\mathbb{Z}_{2^n}[i]$ is $2^{2n-1}-1.$
\end{proposition}
\begin{proof}
   If $R=\mathbb{Z}_{2^n}[i]$, then $Z(R)$ is $\bigcup_{j=1}^{2n-2} V_{d_j}$.
From Lemma $\ref{lem1}(b)$, the number of vertices of $\Gamma(R)$ is $2^{2n-1}-1$.
\end{proof}

\begin{lemma}\label{lem3}
In the zero-divisor graph of $R=\mathbb{Z}_{2^n}[i]$,
\begin{itemize}
\item[$(a)$] $V_{d_j} $ is adjacent to $V_{d_{j'}}$ in $\Gamma(\mathbb{Z}_{2^n}[i])$ if and only if $d_jd_{j'}=0$ in $\mathbb{Z}_{2^n}[i].$
\item[$(b)$] $\displaystyle \bigcup_{j=1}^n V_{d_{j}}$ forms a clique and each $V_{d_{j}}$ is adjacent to $\displaystyle \bigcup_{j'=n+1}^{2n-j}V_{d_{j'}}$, where $j \in[1, n-1]$.
\item [$(c)$] The degree of any vertex $\alpha_k$
    \[ 
\d(\alpha_{k})= \begin{cases} 
      {2^{n+n-j}-2} 
& \text{~if~} \alpha_k \in V_{d_j}, j \in [1,n] \\
      2^{n-j}-1 & \text{~if~}\alpha_k \in V_{d_{n+j}}, j \in [1,n-1],
\end{cases}
\]
hence $\delta(\Gamma(R))=1$, $\Delta(\Gamma(R))=|V(\Gamma(R))|-1$ and the vertex with maximum degree is $2^{n-1}(1+i).$
\end{itemize}
\end{lemma}
\begin{proof}
    \begin{itemize}
\item[$(a)$] Let $u, u'$ be two units in $\mathbb{Z}_{2^n}[i].$ Then $u u' \neq 0.$ Consider two elements $u d_j \in V_{d_j}$, $u' d_{j'} \in V_{d_{j'}}.$ Then $V_{d_j}$ is adjacent to $V_{d_{j'}}$ in $\Gamma(R)$ if and only if $(u d_j)(u' d_j')=0$ if and only if $(uu')(d_jd_{j'})=0$ if and only if $(d_jd_{j'})=0$ in $R$.
\item[$(b)$] The Lemma follows by Lemma $\ref{lem3}(a)$ and definition of $d_j.$
\item[$(c)$] The degree condition of the vertex follows by using Lemma $\ref{lem3}(a)$ and Lemma $\ref{lem1}(a)$. From Lemma $\ref{lem3}(b)$, the set $V_{d_{2n-1}}$ is adjacent only to $V_{d_{1}}$ in $\Gamma(R)$, where $|V_{d_1}|=1$ and hence $\delta(\Gamma(R))=1$. Additionally, $V_{d_{1}}$ is adjacent to $\bigcup_{j=2}^{2n-1}V_{d_{j}}$, where $|V_{d_{1}}|=1$ and $|\bigcup_{j=2}^{2n-1}V_{d_{j}}|=|V(\Gamma(R))|-1.$
    \end{itemize}
\end{proof}

We represent the set of edges of the zero-divisor graph of $R=\mathbb{Z}_{2^n}[i]$ using two subsets $\Lambda$ and $\Omega$ of $\mathcal{Z}(R)$ where $\Lambda=\cup_{j=1}^{n}V_{d_j}$ and $\Omega=\cup_{j=n+1}^{2n-1}V_{d_j}$.
\begin{definition} \label{defn2}
Let $E_{j,j'}$ be the set of all edges $\alpha_j-\alpha_{j'}$, where $\alpha_j \in V_{d_{j}}, \alpha_{j'} \in V_{d_{j'}}$ and $E_{j,j}$ be the set of all edges in each $\langle V_{d_j} \rangle.$ Here note that $E_{j,j'}$ denotes the same set of edges as $E_{j',j}$.\\
The zero-divisor graph of $\mathbb{Z}_{2^n}[i]$ is given by $(V,E)$, where $V=\Lambda \cup \Omega$ and $E=\mathfrak{E}_1 \cup \mathfrak{E}_2$, where $\mathfrak{E_1}$ denotes the set of edges in $\Lambda$  and $\mathfrak{E_2}$ denotes the set of all edges connecting $\Omega$ and $\langle \Lambda \rangle$.
By using the adjacency condition between $V_{d_j}$'s in Lemma $\ref{lem3}(b)$, we get
  \begin{equation*} 
  \begin{split}
\mathfrak{E}_1 =&  \bigcup_{j=1}^{n} \bigcup_{j'=j}^{n} E_{j,j'}\\
    \mathfrak{E}_2 =& \bigcup_{j=1}^{n-1} \bigcup_{j'=1}^{n-j} E_{j,n+j'}.
  \end{split}
 \end{equation*}
 \end{definition}

\begin{theorem}\label{thm1}
The size of the zero-divisor graph of the ring $\mathbb{Z}_{2^n}[i]$ is $2^{2n-1}(n-1)-2^{n-1}+1$.
\end{theorem}
\begin{proof}
Let $R=\mathbb{Z}_{2^n}[i]$. Suppose that the zero-divisor graph $\Gamma(R)$ is $(V,E)$, where the vertex set $V=Z(R)$ and the edge set $E= \mathfrak{E}_1 \cup \mathfrak{E}_2$. \\
Now we proceed to find the number of elements in $\mathfrak{E}_1$ and $\mathfrak{E}_2$. From Lemma $\ref{lem3}(b)$, the number of elements in $\mathfrak{E}_1$ is the size of the complete graph $K_{2^n-1}$ and hence $|\mathfrak{E}_1|=\frac{1}{2}(2^n-1)(2^n-2)=(2^n-1)(2^{n-1}-1)$. \\
By Definition \ref{defn2} and Lemma $\ref{lem1}(a)$,
\begin{equation*} 
\begin{split} 
|\mathfrak{E}_2| = & |\bigcup_{j=1}^{n-1} E_{1,n+j}|+|\bigcup_{j=1}^{n-2} E_{2,n+j}|+\cdots+|\bigcup_{j=1}^{2}  E_{n-2,n+j}|+|\bigcup_{j=1}^{1}E_{n-1,n+j}|\\
  = & n2^{2n-1}-2^{2n}+2^n
\end{split}
\end{equation*}
Therefore, the size of zero-divisor graph of $\mathbb{Z}_{2^n}[i]$ is $2^{2n-1}(n-1)-2^{n-1}+1$.
\end{proof}

\begin{lemma}\label{lem4}
Let $R=\mathbb{Z}_{2^n}[i]$. The clique number of $\Gamma(R)$ is the number of nilpotent elements of index $2$ in $R$.    
\end{lemma}
\begin{proof}
From Lemma $\ref{lem1}(c)$, the clique number is equal to $2^n-1.$
\end{proof}

\begin{proposition}\label{prop3}
The radius and diameter of the zero-divisor graph of $\mathbb{Z}_{2^n}[i]$ are $1$ and $2$, respectively. 
\end{proposition}
\begin{proof}
Let $R=\mathbb{Z}_{2^n}[i]$. Let $\alpha_0=2^{n-1}+2^{n-1}i$. Case 1. If $\alpha_j, \alpha _{j'}\in \Lambda$. Since the graph induced by $\Lambda$ is a complete graph, $d(\alpha_j, \alpha_{j'})=1$. 
    \noindent
    Case 2. If $\beta_j, \beta_{j'} \in \Omega$, then $\beta_j$ is not adjacent to $\beta_{j'}$ in $\Gamma(R)$. However, $\alpha_0$ is adjacent to both $\beta_j$ and $\beta_{j'}$ and therefore $\beta_j-\alpha_0-\beta_{j'}$ is a path in $\Gamma(R)$. Therefore, $d(\beta_j, \beta_{j'})=2$.\\ On the other hand, if $\alpha_j \in  \Lambda$,  and $\beta_{j} \in  \Omega$, then $\alpha_0$ is adjacent to both $\alpha_j$ and $\beta_{j}$ in $\Gamma(R)$. Now, $\alpha_j -\alpha_0-\beta_{j}$ is a path in $\Gamma(R)$ and therefore $d(\alpha_j, \beta_{j})=2$. Thus, the eccentricity of any vertex $\alpha_j$ of $\Gamma(R)$ 
     \begin{equation*}
e(\alpha_j) = 
\left\{
    \begin{array}{lr}
       1  & \text{if~} \alpha_j = \alpha_0 \\
       2 & \text{otherwise}.
    \end{array}
\right.
\end{equation*}
    Hence, the radius and diameter of $\Gamma(R)$ are $1$ and $2$, respectively.
\end{proof}
\begin{corollary}
    $\Gamma(\mathbb{Z}_{2^n}[i])$ is not self-centered, where $n \geq 2$.
\end{corollary}
\begin{proof}
    The vertex $2^{n-1}(1+i)$ has eccentricity 1, and therefore, the subgraph induced by this vertex is $K_1$, which is not isomorphic to $\Gamma(\mathbb{Z}_{2^n}[i])$.
    When $n=1$, $\Gamma(\mathbb{Z}_{2}[i])$ is isomorphic to $K_1$ and is thus self-centered.
\end{proof}

\begin{theorem}\label{thm2}
   For the ring $R=\mathbb{Z}_{2^n}[i]$, the vertex connectivity and edge connectivity of the zero-divisor graph of $R$ are the same and equal to $1$.
\end{theorem}
\begin{proof}
    From Lemma $\ref{lem3}(c)$, removal of the vertex $\alpha_{0}=2^{n-1}(1+i)$ disconnects the graph $\Gamma(R)$ since there are pendant vertices in the graph and $\alpha_{0}$ is adjacent to all other vertices, in particular with pendant vertices. Hence, the vertex connectivity, $\kappa(\Gamma(R))=1.$\\
    Also, from Lemma $\ref{lem3}(d)$, there are vertices in $\Gamma(R)$ having degree 1, that is, each vertex in $V_{d_{2n-1}}$ has degree 1 and in particular $1+i \in V_{d_{2n-1}}$. The removal of the edge, say $e$ having one end vertex as $1+i$ and the other as $2^{n-1}(1+i)$ disconnects the graph, and hence edge connectivity $\lambda(\Gamma(R))=1.$
\end{proof}

\begin{theorem}\label{thm3}
    For the ring $R=\mathbb{Z}_{2^n}[i]$, the zero-divisor graph of $R$ is planar if and only if $n=1$ or $n=2$.
\end{theorem}
\begin{proof}
 From Lemma $\ref{lem3}(b)$, $\cup_{j=1}^{n} V_{d_j}$ forms a clique of order $2^{n}-1$ in $\Gamma(R)$. Suppose that if $n \geq 3$, then the order of the clique is greater than or equal to $7$, and hence $K_5$ is a subgraph of $\Gamma(R)$. By Kuratowski's theorem[\cite{bondymurty}, p.268], $\Gamma(R)$ is not a planar graph.
 Conversely, for $n=1$,  $\Gamma(R)$ is $K_1$, a planar graph, and for $n=2$, $\Gamma(R)$ is a one-point union of $C_3$ and $4K_{2}$, which is again a planar graph. 
\end{proof}


\begin{lemma}\label{lem2.5}
The chromatic number of the subgraph induced by $V_{d_j}$ of $\Gamma(\mathbb{Z}_{2^n}[i])$ is $2^{j-1}$, where $j \in [1,n].$
\end{lemma}
\begin{proof}
    First, we define an ordering for the elements of $V_{d_j}$ by 
    $$V_{d_{j}}=\{\alpha_{j,1}, \alpha_{j,2}, \ldots, \alpha_{j,{2^{j-1}}}: N(\alpha_{j,1}) \leq  N(\alpha_{j,2})  \leq \cdots \leq N(\alpha_{j,2^{j-1}}) \},$$ $V_{d_j}$ is an ordered set with the above ordering. Next, we define the colouring for each element of $V_{d_j}$ through a mapping $f_j: V_{d_j} \longrightarrow [1,2^{j-1}]$ defined by $f_j(\alpha_{j,k})=k.$ Since $f_j$ is a bijection, $f_j$ is a $2^{j-1}$-proper colouring.  By Corollary \ref{cor1} and by equation $14.2$ [\cite{bondymurty}, p.359], the chromatic number of $\langle V_{d_j} \rangle$ is $2^{j-1}$.
    \end{proof}


\begin{theorem}\label{thm4}
 The chromatic number of the zero-divisor graph of $\mathbb{Z}_{2^n}[i]$ is $2^{n}-1$.
\end{theorem}
\begin{proof} Let $R=\mathbb{Z}_{2^n}[i]$. Then $\Gamma(R)=(V,E)$ is a graph with the vertex set $V=\mathcal{Z}(R)$ and $E=\mathfrak{E}_1\cup \mathfrak{E}_2$. Initially, for each element of $V_{d_j}$ for $j \in [1,n]$, the colouring is given as in Lemma \ref{lem2.5}. We proceed to define $f:V(\Gamma(R)) \longrightarrow \{1,2, \ldots, 2^{n} -1\} $ given by
    \[ 
f(\alpha_{j,k})= \begin{cases} 
      1 & \alpha_{1,k} \in V_{d_1} \\
f_j(\alpha_{j,k})+\displaystyle \sum_{\ell=1}^{j-1}|V_{d_{\ell}}|
 & \alpha_{j,k} \in V_{d_j}, j \in [2,n]\\
{2^n-1} 
& \alpha_{j,k} \in  V_{d_j}, j \in [n+1,2n-1]
\end{cases}
\]
We prove it is proper colouring. It is not difficult to verify that $f$ is well-defined. We claim that $f$ is a proper colouring of $\Gamma(R)$, that is, to prove $f(\alpha_k)\ne f(\alpha_{k'})$, whenever $\alpha_k$ is adjacent to $\alpha_{k'}$ in the zero-divisor graph of $R$.\\
Case 1. Suppose $\alpha_{j,k}$ and $\alpha_{j',k'}$ are adjacent, where $\alpha_{j,k} \in V_{d_j}, \alpha_{j',k'} \in V_{d_{j'}}$ for $j,j' \in [1,n]$. Assume $j=j'$. Consider two elements $\alpha_{j,k}, \alpha_{j,k'} \in$  $V_{d_j}$. For fixed $j\in[1,n]$, $\sum_{\ell=1}^{j-1}|V_{d_{\ell}}|$ is a constant and by Lemma \ref{lem2.5}, the set $\{f_j(\alpha_{j,k}): \alpha_{j,k} \in V_{d_{j}}\}$ is a set of consecutive integers, hence  $f_j(\alpha_{j,k})\ne f_j(\alpha_{j,k'})$. On the other hand, if $j<j'$ and $\alpha_{j,k} \in V_{d_j}$, $\alpha_{j',k'} \in V_{d_{j'}}$, then $ \sum_{\ell=1}^{j-1}|V_{d_{\ell}}|<  \sum_{\ell=1}^{j'-1}|V_{d_{\ell}}|$ implies $\max \{f_j(\alpha_{j,k}): \alpha_{j,k} \in V_{d_{j}}\}+  \sum_{\ell=1}^{j-1}|V_{d_{\ell}}| < \min \{f_{j'}(\alpha_{j',k'}): \alpha_{j',k'} \in V_{d_{j'}}\} +  \sum_{\ell=1}^{j'-1}|V_{d_{\ell}}|$ and hence $f(\alpha_{j,k}) \neq f(\alpha_{j',k'})$. Similarly, if $j>j'$ and $\alpha_{j,k} \in V_{d_j}$ and $\alpha_{j',k'} \in V_{d_{j'}}$ then $ \sum_{\ell=1}^{j-1}|V_{d_{\ell}}| >  \sum_{\ell=1}^{j'-1}|V_{d_{\ell}}|$ implies $\min \{f_j(\alpha_{j,k}): \alpha_{j,k} \in V_{d_{j}}\}+ \sum_{\ell=1}^{j-1}|V_{d_{\ell}}| > \max \{f_{j'}(\alpha_{j',k'}): \alpha_{j',k'} \in V_{d_{j'}}\} + \sum_{\ell=1}^{j'-1}|V_{d_{\ell}}|$ and hence $f(\alpha_{j,k}) \neq f(\alpha_{j',k'})$.\\
Case 2. Suppose $\alpha_{j,k}$ and $\alpha_{j',k'}$ are adjacent, where $\alpha_{j,k} \in \{V_{d_{j}}: j\in[1,n] \}$ and $\alpha_{j',k'} \in \{V_{d_{n+j}}: j\in[1,n-1]\}$.
Now, for any $j' \in [1,n-1]$, by the definition of $f$, we have $\max \{f(\alpha_{j',k}): \alpha_{j',k} \in V_{d_{j}}\} < \sum_{\ell=1}^n|V_{d_{\ell}}|=2^n-1$, and hence $f(\alpha_{j,k})\ne f(\alpha_{j',k'})$.
Since $f$ is a proper colouring of $\Gamma(\mathbb{Z}_{2^n}[i])$, $\chi(\Gamma)\leq 2^n-1.$ By Lemma \ref{lem4} and by equation $14.2$ [\cite{bondymurty}, p.359], $\chi(\Gamma) \geq 2^n-1$ and hence $\chi(\Gamma) = 2^n-1.$ This completes the proof.
\end{proof}

\begin{corollary}
The zero-divisor graph of $\mathbb{Z}_{2^n}[i]$ is weakly perfect.
\end{corollary}

\begin{lemma}\label{lem5} 
In the ring $\mathbb{Z}_{2^n}[i]$, there exists a one-to-one correspondence between the set of all $V_{d_j}$'s and $R_E$.
\end{lemma}

\begin{proof}
If $\alpha_1 \in ann(\beta_1)$, then $\alpha_1\beta_1 =0$. Now, $(u\beta_1)\alpha_1=u(\beta_1\alpha_1)=0$ implies $\alpha_1 \in ann(u\beta_1)$. Thus, $ann(\beta_1)=ann(u\beta_1)$ for $u \in \mathcal{U}(R)$. Consider a mapping $\phi:R_E \longrightarrow \{V_{d_j}: 1 \leq j \leq 2n-1\}$ given by $\phi([\beta_j])= V_{d_j}$, where $\beta_j=ud_j$. Suppose $\phi([\beta_j])=\phi([\beta_k])$. Then $V_{d_j}= V_{d_k}$. Suppose $\alpha=u_1d_j \in V_{d_j}=V_{d_k}$ implies $u_1d_j=u_2d_k$, hence $d_j=u_1^{-1}u_2d_k$, it is clear that $d_j$ and $d_k$ are associates. Thus, $[d_j]=[d_k]$ implies $[\beta_j]=[\beta_k].$
\end{proof}
 
\begin{theorem}\label{thm5}
Let $R=\mathbb{Z}_{2^n}[i]$. Then the order and size of $\Gamma_E(R)$ is $2n-1$ and $n^2-n $ respectively.
\end{theorem}
\begin{proof}
From Lemma \ref{lem5}, the order of $\Gamma_E(R)$ is equal to the number of distinct $V_d$'s. It follows from Remark \ref{note1} that the order of $\Gamma_E(R)$ is $2n-1$. The adjacency between $V_d$'s contributes to the set of edges of $\Gamma_E(R)$. From Lemma $\ref{lem3}(b)$,  the size of $\Gamma_E(R)$ equals the sum of size of complete graph of order $n$ and  $n-j$ edges for $j \in [1,n-1]$ implies size of $\Gamma_E(R) = \frac{n(n-1)}{2}+ \displaystyle\sum_{j=1}^{n}(n-j)= n^2-n.$
\end{proof}
     
\section{Matching in the zero-divisor graph of $\mathbb{Z}_{2^n}[i]$}
Further, using structural identification of the zero-divisor graph through associate classes, we determine the independence number, maximum, and maximal matching of the zero-divisor graph of $\mathbb{Z}_{2^n}[i]$ by considering $V_{d_j}$ as an ordered set as defined in Lemma \ref{lem2.5}.

\begin{lemma}{\label{lem3.1}}
   For $j \in [n+1,2n-1]$, each $V_{d_j}$ is an independent set in $\Gamma(\mathbb{Z}_{2^n}[i])$.  Moreover, $\bigcup\limits_{j=n+1}^{2n-1} V_{d_j}$ is an independent set in $\Gamma(\mathbb{Z}_{2^n}[i])$.
\end{lemma}
\begin{proof}
    By Proposition \ref{thm2} and  Lemma $\ref{lem3}(b)$,  $\bigcup\limits_{j=n+1}^{2n-1} V_{d_j}$ is an independent set.
\end{proof}
\begin{lemma}{\label{lem3.2}}
    If $I$ is any maximum independent set of $\Gamma(\mathbb{Z}_{2^n}[i])$, then $I\cap V_{d_n} \ne \emptyset.$
\end{lemma}
\begin{proof}
    $\langle \bigcup\limits_{j=1}^{n} V_{d_j} \rangle$ is a complete subgraph of $\Gamma(\mathbb{Z}_{2^n}[i])$. Hence $I$ can contain only one element from the set $\bigcup\limits_{j=1}^{n} V_{d_j}.$ Using Lemma $\ref{lem3}(b)$ and Lemma \ref{lem3.1}, all elements of $\bigcup\limits_{j=1}^{n-1} V_{d_j}$ are adjacent to a set of independent vertices. To obtain an independent set with maximum cardinality, one must choose that one element from $V_{d_n}$. 
\end{proof}
\begin{lemma}\label{lem3.3}
For each $j\in [2,n]$, let $H_j$ be the subgraph induced by $V_{d_j}$ in $\Gamma(\mathbb{Z}_{2^n}[i])$. Then there exists a perfect matching $M_j$ with $2^{j-2}$ elements in $H_j$.  
\end{lemma}
\begin{proof}
For any $j$, let $C_{2^{j-1}}: \alpha_{j,1}\alpha_{j,2} \cdots \alpha_{j,2^{j-1}-1}\alpha_{j,2^{j-1}}\alpha_{j,1}$ be the largest cycle in $H_j$. Now, consider the following matching $M_j$ in $H_j$.  
   \begin{align*}
       M_j= & \{\alpha_{j,k}\alpha_{j,k+1}: k=2t-1 ~\text{for} ~t \in[1,2^{j-2}]\}
   \end{align*}
   Since graph induced by each $V_{d_j}$ for $j\in [2,n]$ is a complete graph of order $2^{j-1}$ respectively, each $M_j$ is a perfect matching with $\frac{2^{j-1}}{2}=2^{j-2}$ elements.
   \end{proof}
\begin{theorem}\label{thm3.1}
    The vertex independence number of $ \Gamma(\mathbb{Z}_{2^n}[i])$ is $\displaystyle 1+\sum_{j=n+1}^{2n-1} |V_{d_j}|.$
\end{theorem}
\begin{proof}
    Let $I$ be the maximum independent set in $\Gamma(\mathbb{Z}_{2^n}[i])$. The set $I$ is constructed by selecting elements in the order $V_{d_{2n-1}}, V_{d_{2n-2}}, \ldots, V_{d_n}, \ldots, V_{d_2}, V_{d_1}$, since $|V_{d_{2n-1}}| > |V_{d_{2n-2}}| > \cdots > |V_{d_n}| > \cdots > |V_{d_1}|$, in a way such that $I$ is the maximum independent set. From Lemma \ref{lem3.1}, $\bigcup\limits_{j=n+1}^{2n-1}V_{d_j} \in I.$ From Lemma \ref{lem3.2}, an element from $V_{d_n}$ must belong to $I$, say $\alpha_{n,j} \in I$. Since all $V_{d_j}$'s are exhausted, $I= \bigcup\limits_{j=n+1}^{2n-1} V_{d_j} \cup \{\alpha_{n,j}\} \implies |I|= 1+ \sum\limits_{j=n+1}^{2n-1} V_{d_j} .$
\end{proof}

\begin{theorem}\label{thm3.2}
The matching number of $\Gamma(\mathbb{Z}_{2^n}[i])$ is $|V_{d_{n}}|+ |V_{d_{n-1}}| -1$, where $n \geq 2$.  
\end{theorem}
\begin{proof}
For each $j\in [1,n-1]$, consider the following set,
 \begin{align*}
       M_j= & \{\alpha_{j,k}\beta_{j,k}: \alpha_{j,k} \in V_{d_j}, \beta_{j,k} \in V_{d_{n+n-j}}, \textnormal{~and~} k \in [1,2^{j-1}]\}
\end{align*}
with $|M_j|=2^{j-1}$
Since no two edges in $\bigcup\limits_{j=1}^{n-1} M_{j}$ are incident on the same vertex, we get $M_j \cap M_{j'} = \emptyset$ whenever $j \neq j'$ and therefore $\bigcup\limits_{j=1}^{n-1} M_{j}$ is a matching with $\displaystyle \sum_{j=1}^{n-1} |M_{j}|$ number of pairwise disjoint edges. On the other hand, no vertices of the graph $\langle V_{d_{n}} \rangle \cong K_{2^{n-1}} $ is one of the end vertices of any edge in $\bigcup\limits_{j=1}^{n-1} M_{j}$ and therefore, choose 
\begin{align*}
       M_n= & \{\alpha_{n,2t-1}\alpha_{n,2t}: t \in[1,2^{n-2}]\}
   \end{align*}
Consequently,  $M= \bigcup\limits_{j=1}^{n} M_{j}$ is a matching again.\\ 
Since any edge of $\Gamma$ has at least one end vertex in $\Lambda$ and all vertices in $\Lambda$ are saturated by $M$, it is not possible to choose another edge $e$ such that $M \cup \{e\}$ is a matching. Hence, $M$ is a maximum matching with cardinality,
    \begin{equation*}
        \begin{split}
            |M| & =  \displaystyle \sum_{j=1}^{n-1} |M_{j}|+ |M_{n}| =   2^0+ 2^1+\cdots +2^{n-2}+ 2^{n-2}  = 2^{n-1}+ 2^{n-2}-1  =  |V_{d_{n}}|+ |V_{d_{n-1}}|-1.
           \end{split}
    \end{equation*}

 
The following remark provides an alternative proof of the above Theorem using the characterization of maximum matching \cite{bondymurty}, p.416].
\begin{remark}
Consider the same independent edges $M$ as in the proof of Theorem $\ref{thm3.2}$. Now, we claim that this $M$ is the maximum matching. To prove $M$ is maximum, we need to prove there does not exist an augmenting path with respect to this matching $M$; that is, any odd-length alternating path is not an augmenting path with respect to $M.$ Consider an odd length alternating path $P:u-v$ with respect to $M.$ Since the edge, say $e_{11}$ in $M_1$ is in $P$ and by Lemma $\ref{lem3}(d)$, one of the end vertices of $e_{11}$ has degree 1, $e_{11}$ must be incident on either $u$ or $v$.
Assuming without loss of generality that $e_{11}$ is incident to $u$, it follows that $u$ is the matched vertex. Consequently, the alternating path $P$ starts or ends at a matched vertex, specifically at $u$. Hence, $P$ is not an augmenting path.\end{remark}
\end{proof}
\begin{theorem}\label{thm3.4}
    The saturation number of $\Gamma(\mathbb{Z}_{2^n}[i])$ is $|V_{d_n}|-1.$ 
\end{theorem}
\begin{proof}
By Lemma \ref{lem3.3}, let $M_j$ be the set of independent edges from $V_{d_j}$ where $j\in[2,n]$. 

In particular,
\begin{align*}
 M_n= & \{\alpha_{n,2t-1}\alpha_{n,2t}:t \in[1,2^{n-2}]\}
\end{align*}
Consider the edge $e'=\alpha_{1,1}\alpha_{n,1}$. 
Let $M= \bigl(\bigcup_{j=2}^{n} M_j \setminus \{(\alpha_{n,1}\alpha_{n,2})\}\bigr) \cup \{e'\}$. The number of edges in $M$ is given by
\begin{equation*}
|M| = \sum_{j=2}^{n} |M_j|-1 +1 = \sum_{j=2}^{n} 2^{j-2} = 2^{n-1}-1= |V_{d_n}|-1.
\end{equation*}
First, we claim that $M$ is a maximal matching. It is clear that $M$ consists of pairwise nonadjacent edges in $\Gamma(\mathbb{Z}_{2^n}[i])$. To prove $M$ is maximal, it suffices to show that $M \cup \{e\}$ cannot be a matching. Given the construction of $M$, all the vertices in $\Lambda$ are saturated by $M$ except the vertex $\alpha_{n,2}$. For $M \cup \{e\}$ to remain a matching, the end vertices of $e$ must belong to $\Omega$; however, there are no edges connecting the vertices in $\Omega$. Thus, we conclude that $M$ is maximal matching with the cardinality $|V_{d_n}|-1=2^{n-1}-1$. \\
Now, we claim that $s(\Gamma(\mathbb{Z}_{2^n}[i]))=|V_{d_n}|-1$. We know
$\frac{1}{2}(|V(\Gamma)|-\alpha(\Gamma)) \leq s(\Gamma) \leq \alpha'(\Gamma)$. From Theorem \ref{thm3.1} and Theorem \ref{thm3.2}, we have,
\begin{eqnarray*}
     \frac{2^{2n-1}-1-2^{2n-1}+2^n-1}{2} \leq s(\Gamma) \leq 2^{2n-1}+2^{n-2}-1 \\
     2^{n-1}-1 \leq s(\Gamma) \leq 2^{2n-1}+ 2^{n-2}-1
\end{eqnarray*}
\noindent From above inequality, the saturation number of
$\Gamma(\mathbb{Z}_{2^n}[i])$ is $2^{n-1}-1$. 
\end{proof}

\section{Topological indices of the zero-divisor graph of $\mathbb{Z}_{2^n}[i]$}
In this section, we determine important topological indices of the zero-divisor graph of $\mathbb{Z}_{2^n}[i]$.

\begin{theorem}
Let $\Gamma$ be the zero-divisor graph of $\mathbb{Z}_{2^n}[i]$.  Then the Wiener index, the Randi\'c index, the first Zagreb index, and the second Zagreb index of $\Gamma$ are given by 
\begin{itemize}
    \item[$(a)$] $W(\Gamma)=n(2^{2n-1}-2^{2n})+\frac{1}{3}(2^{2n+2}-2^{2n})-(2^{n-1}+2^{2n+1})+2^n+2^{4n-2}+1$ .
    \item[$(b)$] $Rand(\Gamma)=
    \sum_{k=1}^{n-1} \frac{2^{k-1}}{(2^{2n-k}-2)^{1/2}} \bigl( \sum_{j=k+1}^{n} \frac{2^{j-1}}{(2^{2n-j}-2)^{1/2}}+  \sum_{j=1}^{n-k} \frac{2^{(n+j)-1}}{(2^{n-j}-1)^{1/2}}\bigl).$
    \item[$(c)$]
    $M_1(\Gamma)=2^{4n-1}+2^{n+2}+2^{2n-1}-[2^2+n2^{2n+1}+n2^{2n}+2^n].$
    \item[$(d)$]
    $ M_2(\Gamma)= \sum_{k=1}^{n-1} 2^{k-1}(2^{2n-k}-2) \biggl[ \sum_{j=k+1}^{n} (2^{2n-1}-2^j)+ \sum_{j=k+1}^{n} (2^{2n-1}-2^{n+j-(k+1)})\biggl].$
\end{itemize}
\end{theorem}
\begin{proof}
\begin{itemize}
\item[$(a)$] Consider the distance matrix $D_e=(d'_{jj'})$ of the compressed zero-divisor graph $\Gamma_E(\mathbb{Z}_{2^n}[i])$ with the vertices $V_{d_1}, V_{d_2}, \cdots V_{d_{2n-1}}$. Then the distance matrix of the zero-divisor graph of $\mathbb{Z}_{2^n}[i]$ is $D=(d_{jj'})$, where each entry $d_{jj'}$ is a block matrix of order $|V_{d_j}|\times, |V_{d_{j'}}|$ and replace $1$ and $2$ in $d'_{jj'}$ by $\textbf{1}_{|V_{d_j}|\times |V_{d_{j'}}|}$ and $\textbf{2}_{|V_{d_j}|\times |V_{d_{j'}}|}$, respectively whenever $j \neq j'$, and
\begin{equation*}
d_{jj}=\left\{
\begin{array}{ll}
\textbf{1}-\textbf{I}  & \text{for~} j\in [1,n]\\
\textbf{2}- 2\textbf{I}  & \text{for~} j\in [n+1,2n-1],
\end{array}\right.
\end{equation*}
where $\textbf{I}_{|V_{d_j}|\times |V_{d_j}|}$ is the identity matrix and $\textbf{1}_{|V_{d_j}|\times |V_{d_{j'}}|}$ and $\textbf{2}_{|V_{d_j}|\times |V_{d_{j'}}|}$ represent the matrix whose entries are all $1$ and $2$, respectively.
The matrix $D=(d_{jj'})$ is given by,

\begin{align*}
        \begin{array}{c}
            \begin{bNiceMatrix}[first-row,first-col]
            &V_{d_1}& V_{d_2}&\cdots&V_{d_n}&  & V_{d_{n+1}}&\cdots &V_{d_{2n-2}}&V_{d_{2n-1}}\\[6pt]
             V_{d_1}  & \textbf{1-I} & \textbf{1} & \cdots & \textbf{1} &\big| & \textbf{1} & \cdots & \textbf{1} & \textbf{1} \\
V_{d_2} & \textbf{1} & \textbf{1-I} & \cdots & \textbf{1} & \big| & \textbf{1} & \cdots & \textbf{1} & \textbf{2} \\
\vdots& \vdots &\vdots &\ddots & \vdots& \bigg|& \vdots& \cdots & \vdots& \vdots  \\
V_{d_n} & \textbf{1} & \textbf{1} & \cdots & \textbf{1-I} & \big| & \textbf{2} & \cdots & \textbf{2} & \textbf{2} \\
\hline
V_{d_{n+1}}& \textbf{1} & \textbf{1} & \cdots & \textbf{2} & \bigg| &
 \textbf{2}-2\textbf{I} & \cdots & \textbf{2} & \textbf{2}\\
 \vdots & \vdots &\vdots &\cdots & \vdots & \bigg| & \vdots  & \ddots & \vdots& \vdots \\
V_{d_{2n-2}}& \textbf{1} & \textbf{1} & \cdots & \textbf{2} & \big|& \textbf{2}
 & \cdots & \textbf{2}-2\textbf{I} & \textbf{2} \\
V_{d_{2n-1}} &\textbf{1} & \textbf{2} & \cdots & \textbf{2} & \big| & \textbf{2} & \cdots & \textbf{2} & \textbf{2}-2\textbf{I} \\
\end{bNiceMatrix} 
\end{array}
\end{align*}

Since all the entries in $D$ are block matrices, the Wiener index of $\Gamma$ is given by
\begin{equation*}
    \begin{split}
W(\Gamma)= & \frac{1}{2}\sum_{d_{jj'} \in D} \sum_{x \in d_{jj'}}x =  \sum_{\substack{d_{jj'} \in D \\ j < j'}} \sum_{x \in d_{jj'}}x+ \frac{1}{2} \sum_{\substack{d_{jj'} \in D \\ j = j'}} \sum_{x \in d_{jj'}}x =  W_1+ W_2
\end{split}
\end{equation*}
The first term on the right-hand side represents the sum of the elements of each block matrix above the diagonal in matrix $D.$ For $j=j',$ each block matrix is a square matrix, which means the elements below the diagonal in these square matrices lie below the diagonal in $D$, and hence, we take half of its sum.  
Now consider, 
\begin{align*}
    W_1 = & \big[ 1(|V_{d_1}| \times |V_{d_2}|)+ 1(|V_{d_1}| \times |V_{d_3}|)+ \cdots + 1(|V_{d_1}| \times |V_{d_{2n-1}}|) \big]\\
     & + \big[ 1(|V_{d_2}| \times |V_{d_3}|)+ 1(|V_{d_2}| \times |V_{d_4}|)+ \cdots + 1(|V_{d_2}| \times |V_{d_{2n-2}}|) \big] + \cdots + \\
     & + \big[ 1(|V_{d_{n-1}}| \times |V_{d_n}|)+ 1(|V_{d_{n-1}}| \times |V_{d_{n+1}}|) \big]\\
     & + \big[ 2(|V_{d_1}| \times |V_{d_{2n-1}}|) \big]\\
     & + \big[ 2(|V_{d_{2}}| \times |V_{d_{2n-1}}|)+ 2(|V_{d_2}| \times |V_{d_{2n-2}}|) \big]+ \cdots +\\
     & + \big[ 2(|V_{d_{n}}| \times |V_{d_{2n-1}}|)+ 2(|V_{d_n}| \times |V_{d_{2n-2}}|) + \cdots+ 2(|V_{d_n}| \times |V_{d_{n+1}}|)  \big]\\
     & + \big[ 2(|V_{d_{n+1}}| \times |V_{d_{n+2}}|)+ 2(|V_{d_{n+1}}| \times |V_{d_{n+3}}|) + \cdots+ 2(|V_{d_{n+1}}| \times |V_{d_{2n-1}}|)  \big]\\
     & + \big[ 2(|V_{d_{n+2}}| \times |V_{d_{n+3}}|)+ 2(|V_{d_{n+2}}| \times |V_{d_{n+4}}|) + \cdots+ 2(|V_{d_{n+2}}| \times |V_{d_{2n-1}}|)  \big] + \cdots +\\
     & + 2(|V_{d_{2n-2}}| \times |V_{d_{2n-1}}|)\\
     = & \sum_{j=0}^{n-2} \sum_{j'=j+1}^{2n-2-j}  2^{j+j'}+  \sum_{j=1}^{n-1} \sum_{j'=2n-1-j}^{2n-2} 2^{j+j'+1}+ \sum_{j=n}^{2n-3} \sum_{j'=j+1}^{2n-2} 2^{j+j'+1}
\end{align*}
The first term on the right-hand side corresponds to the sum of $1$'s above the main diagonal of $D$, while the second and third terms together represent the sum of the number of $2$'s above the main diagonal of $D$.
Therefore, we have
\begin{align*}
W_1 = &  \sum_{j=0}^{n-2} \sum_{j'=j+1}^{2n-2-j} 2^{j+j'}+  \sum_{j=1}^{n-1} \sum_{j'=2n-1-j}^{2n-2} 2^{j+j'+1}+ \sum_{j=n}^{2n-3} \sum_{j'=j+1}^{2n-2} 2^{j+j'+1}\\
= & 1/3(2+2^{4n-1}-2^{2n+1})+n(2^{2n-1}-2^{2n})-2^{2n}.
\end{align*}
Now we proceed to determine $W_2$.
\begin{align*}
W_2 = & \frac{1}{2} \Bigl( 1  (|V_{d_2}| \times |V_{d_2}|) + 1  (|V_{d_3}| \times |V_{d_3}|)+ \cdots +  (|V_{d_n}| \times |V_{d_n}|) \\
& + 2  (|V_{d_{n+1}}| \times |V_{d_{n+1}}|) + 2  (|V_{d_{n+2}}| \times |V_{d_{n+2}}|)+ \cdots + 2  (|V_{d_{2n-1}}| \times |V_{d_{2n-1}}|)\\
& -  |V_{d_2}| - |V_{d_3}| - \cdots -1  |V_{d_n}|-2 |V_{d_{n+1}}| - 2 |V_{d_{n+2}}| - \cdots - 2 |V_{d_{2n-1}}| \Bigr) \\
= & \frac{1}{3}(2^{2n-1}-2)+ \frac{1}{3}(2^{4n-2}-2^{2n})- 2^{n-1}+1-2^{2n-1}+2^n
\end{align*}
Hence, $W(\Gamma)=n(2^{2n-1}-2^{2n})+\frac{1}{3}(2^{2n+2}-2^{2n})+(2^n-2^{n-1})+(2^{4n-2}-2^{2n+1})+1$.

\item[$(b)$] Consider \begin{equation*} \begin{split}
 Rand(\Gamma) & =  \sum_{\alpha\beta \in E(\Gamma)} \frac{1}{(\text{d}(\alpha)\text{d}(\beta))^{1/2}}\\ & = \sum_{\alpha\beta \in \Lambda} \frac{1}{(\text{d}(\alpha)\text{d}(\beta))^{1/2}}+ \sum_{\alpha\in \Lambda,\beta \in \Omega} \frac{1}{(\text{d}(\alpha)\text{d}(\beta))^{1/2}}\\
 & = Rand_1(\Gamma)+ Rand_2(\Gamma)
\end{split}
\end{equation*}


 Here, the summation in $Rand_1(\Gamma)$ runs over the set of all edges of the graph $\langle \Lambda \rangle$, a complete graph on $2^n-1$ vertices. 
 From Lemma $\ref{lem3}(b)$,
\begin{align*}
  Rand_1(\Gamma) = & \sum_{\substack {\alpha \in V_{d_1},\\ \beta \in V_{d_2}}} \frac{1}{({\text{d}(\alpha) \text{d}(\beta))^{1/2}}}+ \sum_{\substack{\alpha \in V_{d_1}, \\ \beta \in V_{d_3}}} \frac{1}{({\text{d}(\alpha) \text{d}(\beta))^{1/2}}}+ \cdots + \sum_{\substack{\alpha \in V_{d_1}, \\ \beta \in V_{d_n}}} \frac{1}{({\text{d}(\alpha) \text{d}(\beta)})^{1/2}}\\
  & + \sum_{\substack{\alpha \in V_{d_2}, \\ \beta \in V_{d_3}}} \frac{1}{({\text{d}(\alpha)\text{d}(\beta)})^{1/2}}+ \sum_{\substack{\alpha \in V_{d_2}, \\ \beta \in V_{d_4}}} \frac{1}{({\text{d}(\alpha) \text{d}(\beta)})^{1/2}}+ \cdots + \sum_{\substack{\alpha \in V_{d_2}, \\ \beta \in V_{d_n}}} \frac{1}{({\text{d}(\alpha)\text{d}(\beta)})^{1/2}}+ \cdots\\
  & + \sum_{\substack{\alpha \in V_{d_{n-2}}, \\ \beta \in V_{d_{n-1}}}} \frac{1}{(\text{d}(\alpha)\text{d}(\beta))^{1/2}}+ \sum_{\substack{\alpha \in V_{d_{n-2}}, \\ \beta \in V_{d_n}}} \frac{1}{(\text{d}(\alpha) \text{d}(\beta))^{1/2}}+ \sum_{\substack{\alpha \in V_{d_{n-1}}, \\ \beta \in V_{d_n}}} \frac{1}{(\text{d}(\alpha)\text{d}(\beta))^{1/2}}
  \end{align*}
Using Lemma $\ref{lem3}(c)$,
\begin{align*}
Rand_1(\Gamma) = & \frac{|V_{d_1}| |V_{d_2}|}{((2^{2n-1}-2)(2^{2n-2}-2))^{1/2}} + \frac{|V_{d_1}| |V_{d_3}|}{((2^{2n-1}-2)(2^{2n-3}-2))^{1/2}
} + \cdots\\
& + \frac{|V_{d_1}| |V_{d_n}|}{((2^{2n-1}-2)(2^{2n-n}-2))^{1/2}} + \frac{|V_{d_2}| |V_{d_3}|}{((2^{2n-2}-2)(2^{2n-3}-2))^{1/2}} +\\
& + \frac{|V_{d_2}| |V_{d_4}|}{((2^{2n-2}-2)(2^{2n-4}-2))^{1/2}} + \cdots + \frac{|V_{d_2}| |V_{d_n}|}{(2^{2n-2}-2)(2^{2n-n}-2))^{1/2}}+ \cdots \\
& + \frac{|V_{d_{n-2}}| |V_{d_{n-1}}|}{((2^{2n-n+2}-2)(2^{2n-n+1}-2))^{1/2}} + \frac{|V_{d_{n-2}}| |V_{d_n}|}{((2^{2n-n+2}-2)(2^{2n-n}-2))^{1/2}}\\
& + \frac{|V_{d_{n-1}}| |V_{d_n}|}{((2^{2n-n+1}-2)(2^{2n-n}-2))^{1/2}}\\
= & \sum_{{j'}=1}^{n-1} \frac{2^{{j'}-1}}{(2^{2n-{j'}}-2)^{1/2}}\bigg[\sum_{j={j'}+1}^{n}\frac{2^{j-1}}{(2^{2n-j}-2)^{1/2}} \bigg]
\end{align*}
Similarly, using Lemma $\ref{lem3}(b)$ and Lemma $ \ref{lem3}(c)$, we have \\
$Rand_2(\Gamma) =  \displaystyle \sum_{{j'}=1}^{n-1} \frac{2^{{j'}-1}}{(2^{2n-{j'}}-2)^{1/2}}\bigg[\sum_{j=1}^{n-{j'}}\frac{2^{n+j-1}}{(2^{n-j}-1)^{1/2}} \bigg]$\\
Hence, $Rand(\Gamma)=\displaystyle \sum_{{j'}=1}^{n-1} \frac{2^{{j'}-1}}{(2^{2n-{j'}}-2)^{1/2}}\bigg[\sum_{j={j'}+1}^{n}\frac{2^{j-1}}{(2^{2n-j}-2)^{1/2}}+\sum_{j=1}^{n-{j'}}\frac{2^{n+j-1}}{(2^{n-j}-1)^{1/2}} \bigg] $
\item[$(c)$]
    Consider
\begin{equation*}
    \begin{split}
    M_1(\Gamma)= & \sum_{\alpha \in \mathcal{Z}(\mathbb{Z}_{2^n}[i])} \text{d}(\alpha)^2\\
    = & \sum_{\alpha \in \Lambda} \text{d}(\alpha)^2+ \sum_{\alpha \in \Omega} \text{d}(\alpha)^2\\
    \end{split}
    \end{equation*}
    From Lemma $\ref{lem3}(b)$ and Lemma $\ref{lem3}(c)$,
\begin{align*}
    M_1(\Gamma) = & |V_{d_1}| (2^{2n-1}-2)^2+ |V_{d_2}| (2^{2n-2}-2)^2+ \cdots + |V_{d_n}| (2^{2n-n}-2)^2\\
    & + |V_{d_{n+1}}| (2^{n-1}-1)^2+ |V_{d_{n+2}}| (2^{n-2}-1)^2+ \cdots + |V_{d_{2n-1}}| (2^{n-n+1}-1)^2\\
             = & 2^{4n-1}+2^{n+2}+2^{2n-1}-[2^2+n2^{2n+1}+n2^{2n}+2^n].
        \end{align*}
\item[$(d)$]Consider
    \begin{equation*}
    \begin{split}
 M_2(\Gamma) & =  \sum_{\alpha\beta \in E(\Gamma)} \text{d}(\alpha)\text{d}(\beta) \\ & = \sum_{\alpha\beta \in \Lambda} \text{d}(\alpha)\text{d}(\beta)+ \sum_{\alpha\in \Lambda, \beta \in \Omega} \text{d}(\alpha)\text{d}(\beta)\\
 & = Sum_1(\Gamma)+ Sum_2(\Gamma)
 \end{split}
 \end{equation*}
From Lemma \ref{lem3}$(b)$ and Lemma \ref{lem3}$(c)$,  
\begin{align*}
        Sum_1(\Gamma)= & \sum_{\substack{\alpha \in V_{d_1}, \\ \beta \in V_{d_2}}} \text{d}(\alpha) \text{d}(\beta) + \sum_{\substack{\alpha \in V_{d_1}, \\ \beta \in V_{d_3}}} \text{d}(\alpha) \text{d}(\beta)
        + \cdots+ \sum_{\substack{\alpha \in V_{d_1}, \\ \beta \in V_{d_n}}} \text{d}(\alpha) \text{d}(\beta)\\
         & + \sum_{\substack{\alpha \in V_{d_2}, \\ \beta \in V_{d_3}}} \text{d}(\alpha) \text{d}(\beta) + \sum_{\substack{\alpha \in V_{d_2}, \\ \beta \in V_{d_4}}} \text{d}(\alpha) \text{d}(\beta) +  \cdots+ \sum_{\substack{\alpha \in V_{d_2}, \\ \beta \in V_{d_n}}} \text{d}(\alpha) \text{d}(\beta)\\
        & + \sum_{\substack{\alpha \in V_{d_{n-2}}, \\ \beta \in V_{d_{n-1}}}} \text{d}(\alpha) \text{d}(\beta) + \sum_{\substack{\alpha \in V_{d_{n-2}}, \\ \beta \in V_{d_n}}} \text{d}(\alpha) \text{d}(\beta)\\
         & +  \sum_{\substack{\alpha \in V_{d_{n-1}}, \\ \beta \in V_{d_n}}} \text{d}(\alpha)\text{d}(\beta)\\
         = & |V_{d_1}| |V_{d_2}| (2^{2n-1}-2)(2^{2n-2}-2) +|V_{d_1}| |V_{d_3}| (2^{2n-1}-2)(2^{2n-3}-2)\\
         & + \cdots + |V_{d_1}| |V_{d_n}| (2^{2n-1}-2)(2^{2n-n}-2)\\
          & +  |V_{d_2}| |V_{d_3}| (2^{2n-2}-2)(2^{2n-3}-2) +|V_{d_2}| |V_{d_4}| (2^{2n-2}-2)(2^{2n-4}-2)+ \\
         & + \cdots + |V_{d_2}| |V_{d_n}| (2^{2n-1}-2)(2^{2n-n}-2)\\
         & + |V_{d_{n-2}}| |V_{d_{n-1}}| (2^{2n-n+2}-2)(2^{2n-n+1}-2)\\ & +|V_{d_{n-2}}| |V_{d_n}| (2^{2n-n+2}-2)(2^{2n-n}-2)\\
     & + |V_{d_{n-1}}| |V_{d_n}| (2^{2n-n+1}-2)(2^{2n-n}-2)\\
 Sum_1(\Gamma)= & \sum_{k=1}^{n-1}2^{k-1}(2^{2n-k}-2)\bigg[
       \sum_{j=k+1}^{n}(2^{2n-1}-2^j) \bigg]
\end{align*}
Similarly, using Lemma \ref{lem3}$(b)$ and Lemma \ref{lem3}$(c)$, $$Sum_2(\Gamma)= \displaystyle\sum_{k=1}^{n-1}2^{k-1}(2^{2n-k}-2)\bigg[
\sum_{j=k+1}^{n}(2^{2n-1}-2^{n+j-(k+1)}) \bigg]$$\\
Hence, $M_2(\Gamma)= \displaystyle\sum_{k=1}^{n-1} 2^{k-1}(2^{2n-k}-2) \biggl[ \displaystyle\sum_{j=k+1}^{n} (2^{2n-1}-2^j)+ \displaystyle\sum_{j=k+1}^{n} (2^{2n-1}-2^{n+j-(k+1)})\biggl].$
\end{itemize}
\end{proof}
\section{Conclusion}
In this article, we emphasize to identify the structure of the zero-divisor graph of the ring of Gaussian integers modulo $2^n$ via its associate classes of divisors and to determine the chromatic number, maximal and maximum matching. In addition, we obtain a few topological indices of the corresponding zero-divisor graph.  

\section*{Declaration of competing interest}
The authors declare that they have no known competing financial interests or personal relationships that could have
appeared to influence the work reported in this paper.


\section*{Data availability}
No data was used for the research described in the article.
\nocite{*}
\bibliographystyle{plain}
\bibliography{com.bib}
\end{document}